\documentclass[12pt]{article}

\usepackage{amsmath}

\def \R{I\!\!R}
\def \Q{I\!\! \!\!Q}
\def \N{I\!\!N}

\def \esssup{\rm{esssup}}
\def \essmax{\rm{essmax}}
\def \essinf{\rm{essinf}}
\def\INTER{\mathop{\rm {\cap}}\limits}
\newtheorem{definition}{Definition}
\numberwithin{definition}{section}
\newtheorem{theorem}{Theorem}
\newtheorem{lemma}[theorem]{Lemma}
\newtheorem{corollary}[theorem]{Corollary}
\newtheorem{proposition}[theorem]{Proposition}
\newtheorem{remark}[theorem]{Remark}
\newtheorem{notation}[theorem]{Notation}
\numberwithin{theorem}{section}

\begin{document}

\title{Time Consistent Dynamic Risk Processes  \\Cadlag Modification }

\author{Jocelyne Bion-Nadal \\CMAP, Ecole Polytechnique, F-91128 Palaiseau cedex  \\tel: 33 1 69 33 46 25   \\bionnada@cmapx.polytechnique.fr}
\date{July 7, 2006 }

\maketitle

\section*{Abstract}
Working in a continuous time setting, we extend to the general case of dynamic risk measures continuous from above the characterization of time consistency in terms of ``cocycle condition'' of the minimal penalty function. We prove also the  supermartingale property for general time consistent dynamic risk measures.  When the time consistent dynamic risk measure (continuous from above) is normalized and non degenerate, we prove, under a mild condition, that the dynamic risk process of any financial instrument has a cadlag modification. This condition is always satisfied in case of continuity from below.

\section*{Introduction} 
\label{sec:intro}  
This work is in the continuity of our precedent work on dynamic risk measures~\cite{BN02}.
In order to quantify the risk associated to financial positions,  Artzner et al~\cite{Artzner} introduced the notion of  coherent risk measures. This notion was  developed by Delbaen~\cite{D0}. It has been extended  to the convex case by F\"ollmer and Schied~\cite{FS02},~\cite{FS04} and  Frittelli and Rosaza Gianin~\cite{F}, and then to a conditional setting  Detlefsen and Scandolo~\cite{DS}and  Bion-Nadal~\cite{BN01}.
 Coherent dynamic risk measuring has been developed by Delbaen~\cite{D} and Artzner~\cite{ADEHK}. 
Convex dynamic risk measures have been the subject of recent works, see Riedel~\cite{R}, Frittelli and Rosazza Gianin~\cite{FG},  Kl\"oppel and Schweizer~\cite{KS}, Cheredito et al~\cite{CDK}, Bion-Nadal~\cite{BN02} and  F\"ollmer and Penner~\cite{FP}. 
The special case of $g$-expectations or Backward Stochastic Differential Equations has been considered by Peng~\cite{Peng}, Rosazza Gianin~\cite{RG} and Barrieu and El Karoui~\cite{BEK}.

A key property in a dynamic setting is the notion of time consistency. 
It was studied in a continuous time setting for coherent dynamic risk measures by Delbaen~\cite{D}. 
In the convex case, in a discrete time setting, the time consistency  
was characterized by Cheridito et al~\cite{CDK} 
by a condition on the acceptance set and also a ``concatenation condition''.
In a discrete time setting , Detlefsen and Scandolo  ~\cite{DS} have characterized the time consistency for normalized dynamic risk measures by a condition involving an $\essinf$ over all the minimal penalties associated to a family of probability measures.\\  
The characterization in terms of a ``cocycle condition `` 
on the minimal penalty function was first proved by Bion-Nadal~\cite{BN02} 
in a continuous time setting for dynamic convex risk measures continuous  from below. 
It was then proved by F\"ollmer and Penner \cite{FP} for convex normalized 
dynamic risk measures continuous from above, under a restrictive condition 
- namely, that the dynamic risk measure admits a representation in terms of 
probability measures all equivalent to the reference probability $P$. 
Another characterization in terms of a 
supermartingale property can be found, under the same restrictive condition, in
 \cite{FP}. Different notions of time consistency are also studied by Roorda and Schumacher~\cite{RS}.
 
The general setting of the present paper is that of continuous time  and stopping times.
 We define a dynamic risk measure as a family $(\rho_{\sigma,\tau})_{\sigma \leq \tau}$ of convex risk measures on 
$L^{\infty}(\Omega,{\cal F}_{\tau},P)$ 
conditional to $L^{\infty}(\Omega,{\cal F}_{\sigma},P)$ (where $\sigma \leq \tau$ are stopping times).
The time consistency is the composition rule:  $\rho_{\nu,\tau}=\rho_{\nu,\sigma}(-\rho_{\sigma,\tau})$, 
which means that the risk at time $\tau$ associated to a financial position which is represented by 
an essentially ${\cal F}_{\tau}$-measurable function $X$ can be computed either directly 
 or in two steps. 
First we extend (in Section \ref{sec:consistency}) the characterization of time consistency 
in terms of ``cocycle condition'' of the minimal penalty function to the general case of  dynamic risk measures continuous from above. 
The key ingredient in the proof is the use of probability measures $Q$ absolutely continuous with respect 
to the reference probability $P$ and the study $Q$ a.s. of the conditional risk measures $\rho_{\sigma,\tau}$. 
We also extend to the general case, the supermartingale property for the process sum of the dynamic risk process and of the minimal penalty, 
which was introduced in~\cite{FP}.   

Next, Section \ref{sec:cadlag}  
deals with normalized time consistent dynamic risk measure. We prove that if there is  a probability measure $Q$ equivalent to $P$ 
with minimal penalty $0$, the dynamic risk process  associated to each financial 
instrument $\rho_{\sigma,\tau}(X)$  has a cadlag modification which is a $Q$-supermartingale. 
This condition is always satisfied if the dynamic risk measure is non degenerate and continuous from below. 
This result generalizes the one proved by Delbaen in~\cite{D} for coherent risk measures.
In the last Section we give examples of time consistent dynamic risk measures,
using BMO martingales.

\section{ Time consistency}
\label{sec:consistency} 

\subsection{General framework and recalls}

Throughout this paper we work with a filtered probability space\\ 
$(\Omega,{\cal F}_{\infty},({\cal F}_t))_{ t \in \R^+}, P)$. The 
filtration ${\cal F}_t$ is right continuous and ${\cal F}_0$ 
is assumed to be the $\sigma$-algebra generated by the $P$ null sets of ${\cal F}_{\infty}$ 
so that $L^{\infty}(\Omega, {\cal F}_0,P)=\R$. Stopping times are very important in finance. Indeed instant  times can be   defined by  the realization of a particular event ; therefore they are  not deterministic and are modeled as  stopping times.
For every stopping time $\tau$, we consider the $\sigma$-algebra ${\cal F}_{\tau}$ defined by 
${\cal F}_{\tau}=\{A \in \ {\cal F}_{\infty}|\forall t \in \R^+\; A \INTER \{\tau \leq t \} \in {\cal F}_t\}$. 
Let $L^{\infty}(\Omega,{\cal F}_{\tau},P)$ the Banach algebra of essentially bounded  real valued ${\cal F}_{\tau}$-measurable functions. We will always identify an essentially bounded ${\cal F}_{\tau}$-measurable function with its class in $L^{\infty}(\Omega,{\cal F}_{\tau},P)$. A financial position at a stopping time $\tau$ is an element of $L^{\infty}({\cal F}_{\tau})$.\\
We can assume that the time horizon is infinite. Indeed the case of a finite horizon can be considered as a particular case of infinite horizon with ${\cal F}_t={\cal F}_T \;\;\forall t \geq T$.

Recall the following definition of dynamic risk measure (cf~\cite {BN02})
close to the definition of non linear expectations of Peng~\cite{Peng}:  
 
\begin{definition}
$\;$
\begin{enumerate} 
\item A  dynamic risk measure   $(\rho_{\sigma, \tau})_{0 \leq \sigma \leq \tau}$ 
on $(\Omega,{\cal F}_{\infty}, ({\cal F}_t)_{t \in \R^+}, P)$ (where $\sigma \leq \tau$ are two stopping times) is a 
family of maps $(\rho_{\sigma, \tau})_{0 \leq \sigma \leq \tau}$, 
defined on $L^{\infty}({\cal F}_{\tau})$ with values into $L^{\infty}({\cal F}_{\sigma})$ such that each $\rho_{\sigma,\tau}$ is a convex conditional risk measure, i.e. $\rho_{\sigma,\tau}$ satisfies  the following  properties: 
\begin{itemize}
\item[i) ] monotonicity: $$\forall (X,Y) \in (L^{\infty}({\cal F}_{\tau}))^2, \;\;\;if X\leq Y \;\;then \;\;
\rho_{\sigma,\tau}(X) \geq \rho_{\sigma,\tau}(Y)$$  

\item[ii) ] translation invariance: $$\forall Z \in L^{\infty}({\cal F}_{\sigma})\;,\;\;\forall X \in L^{\infty}({\cal F}_{\tau}) 
\;\;\rho_{\sigma,\tau}(X+Z)=\rho_{\sigma,\tau}(X)-Z$$

\item[iii) ] convexity: $$\forall (X,Y) \in (L^{\infty}({\cal F}_{\tau}))^2 \;,\;\;\forall \lambda \in [0,1]$$
$$\rho_{\sigma,\tau}(\lambda X +(1-\lambda)Y) \leq \lambda \rho_{\sigma,\tau}(X)+(1-\lambda) \rho_{\sigma,\tau}(Y)$$
\end{itemize}
\item A dynamic risk measure can have additional properties
\begin{itemize}
\item[iv)]it is normalized if $\rho_{\sigma,\tau}(0)=0\;\; \forall \sigma \leq \tau$
\item[v)]it is continuous from below (resp  above) 
if for every increasing (resp decreasing) sequence $X_n$ of elements of  $L^{\infty}({\cal F}_{\tau})$ such that $X=\lim\;X_n$, the decreasing (resp increasing) 
sequence $\rho_{\sigma,\tau}(X_n) $ has the limit $\rho_{\sigma,\tau}(X)$.
\end{itemize}
\end{enumerate}
\label{definition1}
\end{definition}

\begin{remark}
As it is proved in~\cite{KS} Section 3, and first pointed out by K\"upper, the monotonicity and translation invariance property imply the following regularity property which is usually written in the definition of conditional risk measures: 
$$\forall (X, Y) \in L^{\infty}({\cal F}_{\tau})^2 \;,\; \forall A \in{\cal F}_{\sigma}\;\;
\rho_{\sigma,\tau}( X1_A +Y1_{A^c})=1_A\rho_{\sigma,\tau}( X)+1_{A^c}\rho_{\sigma,\tau}( Y)$$
\label{R1}
\end{remark}

The continuity from below implies continuity from above (cf~\cite{FS04} and~\cite{DS}).

\begin{definition}
The dynamic risk measure is said to be
 time consistent if: \\
$\forall 0 \leq \nu \leq \sigma \leq \tau \;\; \forall X \in L^{\infty}({\cal F}_{\tau})\;\;\; 
\rho_{\nu,\sigma}(-\rho_{\sigma,\tau}(X))=\rho_{\nu,\tau}(X).$ 
\
\label{definition1'}
\end{definition}

 The time consistency condition means that we can  indifferently compute directly the risk at time $\nu$ of a financial position defined at time $\tau$ or compute it in two steps first at time $\sigma$ and then at time $\nu$. 
 It is also related to the following question: Given a normalized monetary risk measure $\rho$ on $(\Omega,{\cal F},P)$, is it possible to ``factorize'' it in a dynamic risk measure or at least through a given stopping time (i.e. $\rho=\rho(-\rho_{\tau,\infty})$)?   

The aim of this  Section is to extend  the  characterization of time consistency  
$\rho_{\nu,\tau}=\rho_{\nu,\sigma}o(-\rho_{\sigma,\tau})$ (for every $\nu \leq \sigma \leq \tau$)  in terms of the ``cocycle  condition'' 
of the minimal penalty functions to the general case of dynamic risk measures continuous from above (i.e. to every dynamic risk measure which admits a representation).\\
Recall that the cocycle condition of the minimal penalty function is defined by the equality
$$\alpha^m_{\nu, \tau}(Q)=\alpha^m_{\nu, \sigma}(Q)+E_Q(\alpha^m_{\sigma, \tau}(Q)|{\cal F}_{\nu})\; \;Q\;a.s.$$for every $Q$.
 This ``cocycle condition'' appeared  in the paper~\cite{BN02} of Bion-Nadal where it is proved (cf theorem 5) 
that for dynamic risk measures continuous from below, the time consistency condition is equivalent to the 
``cocycle condition'' for the minimal penalty function. 
The characterization  of the time consistency condition in terms of ``cocycle condition'' on the minimal penalty function has then 
 been proved, in a discrete time setting,  by F\"ollmer and Penner~\cite{FP} for normalized dynamic risk measures continuous from above,
assuming that there is at least a probability measure $Q$ equivalent to the given probability measure $P$ such that 
$\;\alpha_{0,\tau}^{m}(Q)$ is finite for every stopping time $\tau$.\\
 
The study of the time consistency property for  dynamic risk measures is based upon   the  dual representation of dynamic risk measures. 
Theorems of representation in terms of probability 
measures and of minimal penalties have been proved assuming conditions of continuity from above~\cite{DS},~\cite{BN02}~\cite{FP}.\\

Recall the following results of representation for conditional risk measures

Let two $\sigma$-algebras ${\cal F}_i \subset {\cal F}_j$. Consider a normalized risk measure $\rho_{i,j}$ on $L^{\infty}(\Omega,{\cal F}_j, P)$ conditional to 
$L^{\infty}(\Omega,{\cal F}_i, P)$.\\

\begin{itemize}
\item Dual representation in case of continuity from above proved in ~\cite{DS}

 Every convex risk measure  on  $L^{\infty}(\Omega,{\cal F}_j,P)$ conditional to $L^{\infty}(\Omega,{\cal F}_i,P)$ continuous from above has a representation of the kind:
\begin{equation}
\forall X \in L^{\infty}({\cal F}_{j}) \;\;\rho_{i,j}(X)= \esssup_{Q \in  \tilde{\cal M}_{i,j}}((E_Q(-X|{\cal F}_{i})-\alpha^m_{i,j}(Q))
\label{eq_2}    
\end{equation}
where $\tilde{\cal M}_{i,j}=\{Q\; on\; (\Omega,{\cal F}_j)\;|Q\ll P\; ,\; Q_{|{\cal F}_i}=P\}$.

\item Dual representation in case of continuity from below (~\cite{BN01})

 Every convex risk measure $\rho_{i,j}$ on $L^{\infty}(\Omega,{\cal F}_j, P)$ conditional to 
$L^{\infty}(\Omega,{\cal F}_i, P)$  continuous from below  admits a representation of the kind 
(cf~\cite{BN02}):
\begin{equation}
\forall X \in L^{\infty}({\cal F}_{j}) \;\;\rho_{i,j}(X)= \essmax_{Q \in {\cal M}_{i,j}}((E_Q(-X|{\cal F}_{\sigma})-\alpha^m_{i,j}(Q))
\label{eq_1}    
\end{equation}
 where ${\cal M}_{i,j}=\{Q\; on \;(\Omega,{\cal F}_j)\;|Q\ll P\; ,\; Q_{|{\cal F}_i}=P\; and\; \alpha^m_{i,j}(Q)\in L^{\infty}(\Omega, {\cal F}_i,P)\}$. 
\end{itemize}
Notice that in this representation of dynamic risk measures continuous from below, there are two very useful points:\\
- the fact that the representation is expressed as an ``~$\essmax$~'' and not only as an ``~$\esssup$~''. Indeed,  for every 
 $X \in L^{\infty}({\cal F}_{j})$, there is $ Q_X \in {\cal M}_{i,j}$
 such that $\rho_{i,j}(X)= E_{Q_X}(-X|{\cal F}_{i})-\alpha^m_{i,j}(Q_X)$;\\ 
- and consequently the fact that it is expressed in term of probability measures $Q$ in ${\cal M}_{i,j}$,
i.e. probability measures for which the penalty $\alpha^m_{i,j}(Q)$ is essentially bounded.\\
These two facts were used in~\cite{BN02} in order to prove, for dynamic risk measures constinuous from below, the characterization of time consistency in terms of the ``cocycle condition'' of the minimal penalty.\\
The proof of the characterization of time consistency in~\cite{FP} is given under the condition that the dynamic  risk measure admits a representation in terms of probability measures all equivalent to $P$. This hypothesis is fundamental in the  proof of  ~\cite{FP}. Furthermore that proof  is done only for normalized risk measures and uses  the non negativity of the penalties.

In order to extend the characterization of time consistency in terms of ``cocycle condition'' to the general case of dynamic risk measures continuous from above, the probability measures absolutely continuous with respect to $P$ will play a crucial role. We have to study the conditional risk measures $\rho_{i,j}$ $Q\;a.s.$ for any probability measure $Q$ absolutely continuous with respect to $P$, and in particular we will have to prove a representation theorem for the projection of $\rho_{i,j}$ onto $L^{\infty}(Q)$.

\subsection{Extension of the theorem of representation} 
Let $\rho_{i,j}$ a risk measure $\rho_{i,j}$ on  $L^{\infty}(\Omega, {\cal F}_j,P)$
conditional to  $L^{\infty}(\Omega, {\cal F}_i,P)$. Given a probability measure $Q$ on $(\Omega,{\cal F}_j)$ absolutely continuous with respect to $P$, we want to obtain a representation if $\rho_{i,j}$ $Q$ a.s. i.e. a representation of the projection of $\rho_{i,j}$ onto  $L^{\infty}(\Omega, {\cal F}_i,Q)$. This representation cannot be deduced from the usual representation theorem for risk measures on $L^{\infty}(\Omega,{\cal F}_j,Q)$ conditional to $L^{\infty}(\Omega,{\cal F}_i,Q)$. Indeed   if two random variables $X$ and $Y$ , $(\Omega, {\cal F}_j)$ measurable are equal $Q$ almost surely,  it is possible that the random variables  $\rho_{i,j}(X)$ and $\rho_{i,j}(Y)$ are not equal $Q$ a.s.\\
Introduce the following notations extending usual ones:
 \begin{notation}
 For every probability measure $Q$ on $(\Omega,{\cal F}_j)$ absolutely continuous with respect to $P$,
denote:
\begin{itemize}
\item[i)] Q acceptance set 
\begin{equation}
{\cal A}_{i,j}(Q)=\{Y \in  L^{\infty}(\Omega, {\cal F}_j,P)\;|\;\rho_{i,j}(Y) \leq 0\;\;Q \;a.s.\}
\label{eq_5}    
\end{equation}
\item[ii)] minimal penalty given (as justified in the next remark) by one of the 3 equivalent formula
\begin{eqnarray}
\alpha^m_{i,j}(Q) & = & Q \; ess \; sup_{X \in L^{\infty}(\Omega, {\cal F}_j,P)} 
(E_Q(-X| {\cal F}_i)-\rho_{i,j}(X)) \nonumber \\ 
 & = & Q \; ess \;sup_{Y \in {\cal A}_{i,j}(Q)}  E_Q(-Y| {\cal F}_i) \\ 
 & = & Q \; ess \;sup_{Y \in {\cal A}_{i,j}}  E_Q(-Y| {\cal F}_i) \nonumber 
\label{eq_6}    
\end{eqnarray}
\item[iii) ] $\tilde{{\cal M}_{i,j}}(Q)=\{R \ll P \;|\;R_{|{\cal F}_{i}}=Q\}$
\item[iv) ] ${\cal M}^1_{i,j}(Q)=\{R \in \tilde{{\cal M}_{i,j}} \;|\; E_Q(\alpha^m_{i,j}(R))=E_R(\alpha^m_{i,j}(R)) < \;\infty\} $
\end{itemize}
\label{notation1}
\end{notation}
\begin{remark}
1) The equalities of the equation (4) follow easily from the inclusions
${\cal A}_{i,j} \subset  {\cal A}_{i,j}(Q) \subset L^{\infty}(\Omega, {\cal F}_j,P)$,
and from the fact that, for every $X \in  L^{\infty}(\Omega, {\cal F}_j,P)$,
$ X+\rho_{i,j}(X) \in {\cal A}_{i,j}$.\\
2) When $Q$ is equal to $P$ we have the usual definition of acceptance set, minimal penalty... and in that case we will omit $P$ in the notation.
\label{rem1}
\end{remark}
In the following Lemma we prove that 
for every probability measure $Q$ absolutely continuous with respect to $P$ the canonical projection of $\rho_{i,j}(X) $ onto $L^{\infty}(\Omega, {\cal F}_i,Q)$ can be represented  in terms of probability measures absolutely continuous with respect to $P$ such that their restriction to ${\cal F}_i$ is equal to $Q$.  
\begin{lemma}
For every probability measure $Q$ absolutely continuous with respect to $P$, 
\begin{itemize}
\item[i) ] \begin{equation}
\rho_{i,j}(X)=Q \; ess \;sup_{R \in {\cal M}_{i,j}^1(Q)}(E_R(-X|{\cal F}_i)-\alpha^m_{i,j}(R)) \;\;Q\;a.s.
\label{eq_7}    
\end{equation}
\item[ii) ] For every $X \in  L^{\infty}(\Omega, {\cal F}_j,P)$, there is a sequence $Q_n$ of probability measures in ${\cal M}_{i,j}^1(Q)$ such that $Q\;a.s.$, $\rho_{i,j}(X)$ is the increasing limit of the sequence   $E_{Q_n}(-X|{\cal F}_i)-\alpha^m_{i,j}(Q_n)$
\item[iii) ] There is a sequence $Z_n \in {\cal A}_{i,j}(Q)$ such that $\alpha^m_{i,j}(Q)$
is $Q$ a.s. the increasing limit of the sequence $E_Q(-Z_n|{\cal F}_i)$.
\end{itemize}
\label{lemma0}
\end{lemma}
Before giving the proof of this lemma we recall the two following results:
\begin{lemma}
\begin{enumerate}
\item Let $Q$ a probability measure. Assume that $K$ is a lattice upward directed. Then $Q \esssup \{X \in K\}$ is the limit $Q$ a.s. of an increasing sequence $X_n$ of elements of $K$
\item Let $(\Omega,{\cal F}_j,P)$ a probabilty space . Let ${\cal F}_i$ a subsigma algebra of ${\cal F}_j$. Let  $Z_n$  an increasing sequence of essentially bounded ${\cal F}_j$-measurable functions converging $Q$ a.s. to $Z$. Assume that $E(Z_1)>-{\infty}$. Then $E(Z_n)$ tends to $E(Z)$ and $E(Z_n|{\cal F}_i)$ converges $Q$ a.s. to $E(Z|{\cal F}_i)$. 
\end{enumerate}
\label{lemmaS}
\end{lemma}
 For the  first part  of the lemma, we refer to the appendix of ~\cite{FS02}\\
 For  the second part we refer to ~\cite{DM1}
 
{\bf Proof} of Lemma \ref{lemma0}.  
 We begin with the proof of iii) that we will use in order to prove i).

{\it Proof} of iii). The set $\{E_Q(-Z|{\cal F}_i)\;|\; Z \in {\cal A}_{i,j}(Q)\}$ is a lattice upward directed. Thus we get the existence of $Z_n$,applying Lemma \ref{lemmaS}  1)    

{\it Proof} of i). We adapt the proofs of representation of conditional risk measures~\cite{DS} ,~\cite{BN01} and~\cite{FP}. To do this part of the proof,  we can assume that $\rho_{i,j}(0)=0$.  Indeed the representation result is satisfied for $\rho_{i,j}$ if and only if it is satisfied for $\rho_{i,j}-\rho_{i,j}(0)$.
 The inequality 
$$\rho_{i,j}(X) \geq Q \; ess \;sup_{R \in {\cal M}_{i,j}^1(Q)}(E_R(-X|{\cal F}_i)-\alpha^m_{i,j}(R)) \;\;Q\;a.s.$$ follows from the definition of $\alpha^m_{i,j}(R)$.
Thus in order to prove (\ref{eq_7}), it is enough to verify that
$$E_Q(\rho_{i,j}(X)) \leq E_Q (\; ess \;sup_{R \in {\cal M}_{i,j}^1(Q)}(E_R(-X|{\cal F}_i)-\alpha^m_{i,j}(R)).$$
Define $\rho(X)=E_Q(\rho_{i,j}(X))$. $\rho$ is a monetary risk measure continuous from above thus from theorem 4.31 of~\cite{FS04}, $$\rho(X)=sup_{\{R \;|\alpha^m(R)<\infty\}}(E_R(-X)-\alpha^m(R))$$

As $\rho_{i,j}$ is normalized, the restriction of $\rho$ to $L^{\infty}({\cal F}_i)$ is equal to $E_Q$ and therefore the restriction of every $R$ to ${\cal F}_i$ has to be equal to $Q$ indeed:\\ let $R$ such that $\alpha^m(R)$ is finite.
As $\forall \beta \in \R\;\;\forall B \in {\cal F}_i$, $\rho_{i,j}(\beta 1_B)=
-\beta 1_B$, we get $\beta (Q(B)-R(B)) \leq \alpha^m(R)\;\;\forall \beta \in \R$. So $Q(B)=R(B)\;\;\forall B \in {\cal F}_i$ i.e. $ R \in \tilde {\cal M}_{i,j}(Q)$.\\
Now we compare $\alpha^m(R)$ with $ E_R(\alpha^m_{i,j}(R))$.\\
As  $\alpha^m(R)=sup_{Y \in L^{\infty}({\cal F}_j)} (E_R(-Y)-\rho(Y))$
it follows from the definition of $\alpha^m_{i,j}(R))$ (equation (4)) that $\alpha^m(R) \leq E_R(\alpha^m_{i,j}(R))$.\\
We apply property iii) (already proved) of the lemma to $R$. From Lemma \ref{lemmaS} 2), it follows that $E_R(\alpha^m_{i,j}(R))$ is the limit of 
$E_R(-Z_n)$. $Z_n$ is in $ {\cal A}_{i,j}(R) \subset {\cal  A}_{\rho}$. And this gives the other inequality  $E_R(\alpha^m_{i,j}(R))\leq \alpha^m(R)$ and thus $R \in {\cal M}^1_{i,j}(Q)$ and i) is proved.

{\it Proof} of ii). $\{(E_{R}(-X|{\cal F}_i)-\alpha^m_{i,j}(R)\;|R \in {\cal M}_{i,j}^1(Q)\}$
is a lattice upward directed. Indeed:  Let $R_1, R_2 \; \in {\cal M}_{i,j}^1(Q)$, and denote\\
$A= \{\omega |(E_{R_1}(-X|{\cal F}_i)-\alpha^m_{i,j}(R_1) \geq (E_{R_2}(-X|{\cal F}_i)-\alpha^m_{i,j}(R_2)\}$.\\
The probability measure $R=R_1 1_A+R_2 1_{A^c}$ is in ${\cal M}_{i,j}^1(Q)$ and satisfies\\
$E_{R}(-X|{\cal F}_i)-\alpha^m_{i,j}(R)=sup (E_{R_1}(-X|{\cal F}_i)-\alpha^m_{i,j}(R_1),E_{R_2}(-X|{\cal F}_i)-\alpha^m_{i,j}(R_2))$.\\
So we get the result applying Lemma \ref{lemmaS} 1).

Now we are able to prove the characterization of the time consistency in terms of the cocycle condition in the general setting of dynamic risk measures continuous from above.

\subsection{Characterization of the time consistency}

\begin{theorem}
Consider $\rho_{\sigma,\tau}$ a dynamic risk measure continuous from above.
Let $\nu \leq \sigma \leq \tau$ be stopping times.
The three following conditions are equivalent:
\begin{itemize} 
\item[i) ] The dynamic risk measure is time consistent i.e.
$$\rho_{\nu,\tau}(X)=\rho_{\nu,\sigma}(-\rho_{\sigma,\tau}(X)) \;\;\forall X \in L^{\infty}(\Omega, {\cal F}_{\tau})$$
\item[ii) ] For every probability measure $Q$ absolutely continuous with respect to~$P$,
$${\cal A}_{\nu,\tau}(Q)={\cal A}_{\nu,\sigma}(Q)+{\cal A}_{\sigma,\tau}.$$
\item[iii) ] For every probability measure $Q$ absolutely continuous with respect to $P$, 
the minimal prenalty function satisfies the following cocycle condition 
$$\alpha^m_{\nu, \tau}(Q)=\alpha^m_{\nu, \sigma}(Q)+E_Q(\alpha^m_{\sigma, \tau}(Q)|{\cal F}_{\nu})\; \;Q\;a.s.$$
\end{itemize}
\label{thm1}
\end{theorem}
The equivalence of the two first properties can be found in ~\cite{CDK}.
We deduce this theorem from the following proposition on composition of normalized conditional risk measures:

Consider three $\sigma$-algebras ${\cal F}_1 \subset {\cal F}_2 \subset
{\cal F}_3 $ on a space $\Omega$. Let $P$ a probability measure on $(\Omega,{\cal F}_3)$. For ${1 \leq i \leq j \leq 3}$, let $(\rho_{i,j})$ a  normalized risk measure on $L^{\infty}(\Omega,{\cal F}_j,P)$  conditional to $L^{\infty}(\Omega,{\cal F}_i,P)$.
\begin{proposition}
Assume that $(\rho_{i,j})_{1 \leq i \leq j \leq 3}$ are continuous from above.
The following properties are equivalent:
\begin{itemize} 
\item[i) ] $\rho_{1,3}(X)=\rho_{1,2}(-\rho_{2,3}(X))\; \;\forall X \in L^{\infty}(\Omega,{\cal F}_3,P)$
\item[ii) ] For every probability measure $Q$ absolutely continuous with respect to~$P$,
$${\cal A}_{1,3}(Q)={\cal A}_{1,2}(Q)+{\cal A}_{2,3}.$$
\item[iii) ] For every probability measure $Q$ absolutely continuous with respect to $P$, the minimal penalty function satisfies the following cocycle condition:
\begin{equation}
\alpha^m_{1,3}(Q)=\alpha^m_{1,2}(Q)+E_Q(\alpha^m_{2,3}(Q)|{\cal F}_1) \;Q\;a.s.
\label{eq_8}    
\end{equation}
\end{itemize}
\label{prop2bis}
\end{proposition}
{\bf Proof}:
We adapt the proof of theorem 2 of~\cite{BN02} which was given in the case of continuity from below.

- i) implies ii): \\
Let $X \in {\cal A}_{1,3}(Q)$,  $\rho_{1,3}(X) \leq 0 \; Q\;a.s. $ 
Denote $Z=X+\rho_{2,3}(X)$. By translation invariance, $\rho_{2,3}(Z)=0$ so $Z \in {\cal A}_{2,3}$\\
$\rho_{1,2}(X-Z)=\rho_{1,2}(-\rho_{2,3}(X))=\rho_{1,3}(X) \leq 0 \;Q \;a.s.$ So $X-Z \in  {\cal A}_{1,2}(Q)$. Hence  ${\cal A}_{1,3}(Q) \subset {\cal A}_{1,2}(Q)+ {\cal A}_{2,3}$.\\
Conversely let $Y \in  {\cal A}_{1,2}(Q),\; Z \in {\cal A}_{2,3}$. 
$\rho_{1,3}(Y+Z)=\rho_{1,2}(-\rho_{2,3}(Z)+Y)$.\\
As $Z \in {\cal A}_{2,3}$, $-\rho_{2,3}(Z)+Y \geq Y \; P\;a.s.$. By monotonicity,\\
$\rho_{1,2}(-\rho_{2,3}(Z)+Y) \leq \rho_{1,2}(Y) \leq 0 \;Q \;a.s.$ 
and hence $Y+Z \in {\cal A}_{1,3}(Q)$.\\
Thus ii) is proved.

- ii) implies iii): \\
From the equation (4), 
$$\alpha^m_{1,3}(Q)=Q \; ess \;sup_{X \in {\cal A}_{1,3}(Q)} (E_Q(-X| {\cal F}_1))$$
From ii), it follows that 
$$\alpha^m_{1,3}(Q)=Q \; ess \;sup_{Y \in {\cal A}_{1,2}(Q)}  E_Q(-Y| {\cal F}_1)+Q \; ess \;sup_{Z \in {\cal A}_{2,3}}  E_Q(-Z| {\cal F}_1)$$
Thus it only remains to prove that 
\begin{equation}
Q \; ess \;sup_{Z \in {\cal A}_{2,3}}  E_Q(-Z| {\cal F}_1)=E_Q(\alpha^m_{2,3}(Q)|{\cal F}_1 )\;Q \;a.s.
\label{eq_9}    
\end{equation}
For every $Z \in  {\cal A}_{2,3}$,  $ E_Q(-Z| {\cal F}_2) \leq \alpha^m_{2,3}(Q)\;Q \;a.s.$ 
So we get the inequality $Q \; ess \;sup_{Z \in {\cal A}_{2,3}}   E_Q(-Z| {\cal F}_1)  \leq  E_Q(\alpha^m_{2,3}(Q)|{\cal F}_1 )\;Q \;a.s.$.\\
 Now  exactly as in the  equation (4), we have the equality 
$$Q \; ess \;sup_{Z \in {\cal A}_{2,3}}  E_Q(-Z| {\cal F}_1)=Q \; ess \;sup_{Z \in {\cal A}_{2,3}(Q)}  E_Q(-Z| {\cal F}_1)$$
 From lemma \ref{lemma0}, there is a sequence 
$Z_n\in {\cal A}_{2,3}(Q)$ such that $\alpha^m_{2,3}(Q)$
is $Q$ a.s. the increasing limit of the  sequence 
$E_Q(-Z_n|{\cal F}_2)$. Applying Lemma \ref{lemmaS} 2),  it follows that\\ 
 $E_Q(\alpha^m_{2,3}(Q)|{\cal F}_1)$ is $Q$ a.s. the limit of $E_Q(-Z_n|{\cal F}_1)$  and thus \\
$E_Q(\alpha^m_{2,3}(Q)|{\cal F}_1) \leq Q \; ess \;sup_{Z \in {\cal A}_{2,3}(Q)}  E_Q(-Z| {\cal F}_1)$\\
and this proves the equation (\ref{eq_9}).

- iii) implies i): \\
 This part of the proof  is close to the proof of theorem 3 of~\cite{BN02}.\\
 From lemma \ref{lemma0}, for $X$ fixed, there is a sequence $R_n \in {\cal M}^1_{1,2}(P)$ such that  $\rho_{1,2}(-\rho_{2,3}(X))$ is the increasing limit of $E_{R_n}(\rho_{2,3}(X)|{\cal F}_1)-\alpha^m_{1,2}(R_n)$.\\
- From lemma \ref{lemma0}, for $n$ fixed, there is then a sequence $(Q^n_k)_{k \in \N}$ of probability measures in ${\cal M}^1_{2,3}(R_n)$ such that 
$\rho_{2,3}(X)$ is $R_n \;a.s.$ the limit of the increasing sequence
$E_{Q^n_k}(-X|{\cal F}_2)-\alpha^m_{2,3}(Q^n_k)$.\\

As $Q^n_1 \in {\cal M}^1_{2,3}(R_n)$, $E_R(\alpha^m_{2,3}(Q^n_1)) < \infty$ and from Lemma \ref{lemmaS}, it follows that \\
$(E_{R_n}(\rho_{2,3}(X)|{\cal F}_1)-\alpha^m_{1,2}(R_n))=lim_{k \rightarrow \infty} (E_{R_n}(E_{Q^n_k}(-X|{\cal F}_2)-\alpha^m_{2,3}(Q^n_k)|{\cal F}_1)- \alpha^m_{1,2}(R_n))\;Q\;a.s.$\\
Applying the hypothesis iii) we get\\ 
$E_{R_n}((\rho_{2,3}(X)|{\cal F}_1)-\alpha^m_{1,2}(R_n))=lim_{k \rightarrow \infty}  (E_{Q^n_k}(-X|{\cal F}_1)-\alpha^m_{1,3}(Q^n_k))\;Q\;a.s.$\\
 It follows that $\rho_{1,2}(-\rho_{2,3}(X)) \leq \rho_{1,3}(X)$.\\
- Conversely from lemma \ref{lemma0}, there is a sequence $Q_n \in {\cal M}^1_{1,3}(P)$ such that $P\;a.s.$ $\rho_{1,3}(X)$ is  $P\;a.s.$ the limit of the increasing sequence $(E_{Q_n}(-X|{\cal F}_1)-\alpha^m_{1,3}(Q_n))$.\\
From hypothesis iii)  we get that\\
 $\rho_{1,3}(X)$ is $P\;a.s.$ the limit of $((E_{Q_n}(E_{Q_n}(-X|{\cal F}_2)-\alpha^m_{2,3}(Q_n)|{\cal F}_1)-\alpha^m_{1,2}(Q_n))$.\\
 This gives the converse inequality
$\rho_{1,3}(X) \leq \rho_{1,2}(-\rho_{2,3}(X))$.\\
q.e.d.

Another  useful and constructive result concerning this ``cocycle condition'' is the following (cf~\cite{BN02}): 
Every  convex dynamic risk measure  constructed from a stable family of equivalent probability measures and a penalty function, 
which is both local and satisfies the cocycle condition  defines a time consistent dynamic risk measure. This result extends to the convex case the  result proved by Delbaen~\cite{D} for coherent dynamic risk measure. This result for convex dynamic risk measures is 
important  because it doesn't assume that the penalty function is the minimal one. It allows for the construction of new families of time consistent dynamic risk measures with possible jumps and generalizing the Backward Stochastic Differential Equations~\cite{BN02}.\\

We extend now to general time consistent dynamic risk measure continuous from above the supermartingale property which was proved in~\cite{FP} in the case where the risk measure is normalized and admits a representation in terms of probability measures all equivalent to the reference probability measure. 
 
\begin{proposition}
Let $\rho_{\sigma,\tau}$ a dynamic risk measure continuous from above. It is time consistent,if and only if  for every stopping times $\nu \leq \sigma \leq \tau$, the following two properties are satisfied
\begin{itemize} 
\item[1) ] supermartingale property:\\
for every probability measure $Q$ absolutely continuous with respect to~$P$, such that $E_Q(\alpha^m_{\nu,\tau}(Q))$ is finite, 
$\forall X \in L^{\infty}(\Omega,{\cal F}_{\tau}, P)$, 
\begin{equation}
E_Q(\rho_{\sigma,\tau}(X)+\alpha^m_{\sigma,\tau}(Q)|{\cal F}_{\nu}) \leq 
\rho_{\nu,\tau}(X)+\alpha^m_{\nu,\tau}(Q)\;\;Q\;a.s.
\label{eq_S}
\end{equation}
\item[2) ]
\begin{equation}
\forall Z \in L^{\infty}({\cal F}_{\sigma})\;\;\rho_{\nu,\sigma}(Z)=\rho_{\nu,\tau}(Z+\rho_{\sigma,\tau}(0))
\label{eq_N}
\end{equation}
\end{itemize}

\end{proposition}

\begin{remark}
When the dynamic risk measure is normalized the condition 2) means that $\rho_{\nu,\sigma}$ is the restriction of $\rho_{\nu,\tau}$ to $L^{\infty}({\cal F}_{\nu})$.
\end{remark}

{\bf Proof}:

- Assume time consistency.
Let $Q$ a probability measure absolutely continuous with respect to $P$.
From lemma \ref{lemma0} ii), there is a sequence $Q_n \in {\cal M}^1_{\sigma,\tau}(Q)$ such that   
$\rho_{\sigma,\tau}(X)$ is $Q\;a.s.$ the increasing limit of $E_{Q_n}(-X|{\cal F}_{\sigma})-\alpha^m_{\sigma,\tau}(Q_n)$. As $E_Q(\alpha^m_{\sigma,\tau}(Q_1))<\infty$, we can apply Lemma \ref{lemmaS} 2) . Thus $E_Q(\rho_{\sigma,\tau}(X)|{\cal F}_{\nu})$ is $Q\;a.s.$ the limit of 
$E_Q(E_{Q_n}(-X|{\cal F}_{\sigma})-\alpha^m_{\sigma,\tau}(Q_n)|{\cal F}_{\nu})$. For every $n$, $\alpha^m_{\sigma,\tau}(Q_n)$ is bounded from below and 
$E_Q(\alpha^m_{\sigma,\tau}(Q_n))$ is finite, thus $\alpha^m_{\sigma,\tau}(Q_n)$ is finite $Q$ almost surely. It is the same for $\alpha^m_{\nu,\tau}(Q)$. Applying the cocycle condition to $Q$ and $Q_n$,  we get the existence of a set $A$ such that $Q(A^c)=0$ and such that $\alpha^m_{\sigma,\tau}(Q_n)$,  $\alpha^m_{\nu,\tau}(Q)$, $\alpha^m_{\nu,\sigma}(Q)$ , $E_Q(\alpha^m_{\sigma,\tau}(Q)| {\cal F}_{\nu})$, $\alpha^m_{\nu,\tau}(Q_n)$ are  finite on $A$. We apply the cocycle condition to $Q_n$, as the restriction of $Q_n$ to ${\cal F}_{\sigma}$ is equal to $Q$, and as every term is finite on $A$ and $Q(A^c)=0$, we get 
$$E_Q(E_{Q_n}(-X|{\cal F}_{\sigma})-\alpha^m_{\sigma,\tau}(Q_n)|{\cal F}_{\nu})
=E_{Q_n}(-X|{\cal F}_{\nu})-\alpha^m_{\nu,\tau}(Q_n)+\alpha^m_{\nu,\sigma}(Q)\;Q\;a.s.$$
And thus 
$E_Q(\rho_{\sigma,\tau}(X)|{\cal F}_{\nu})\leq \rho_{\nu,\tau}(X)+\alpha^m_{\nu,\sigma}(Q)\;Q\;a.s.$\\
Applying now the cocycle condition for $Q$, we get the required inequality (\ref{eq_S}). The equality (\ref{eq_N}) is obvious.\\
- Conversely, assume that inequality (\ref{eq_S}) is satisfied for every probability measure $Q \in {\cal M}^1_{\nu,\tau}$
Let $(Y,Z) \in (L^{\infty}(\Omega,{\cal F}_{\tau})^2$ such that $\rho_{\sigma,\tau}(Y)=\rho_{\sigma,\tau}(Z)$. 
Applying (\ref{eq_S}) we get
\begin{eqnarray}
\rho_{\nu,\tau}(Y)+\alpha^m_{\nu,\tau}(Q) & \geq & E_Q(\rho_{\sigma,\tau}(Y)||{\cal F}_{\nu})+ E_Q(\alpha^m_{\sigma,\tau}(Q)|{\cal F}_{\nu}) \nonumber \\
\rho_{\sigma,\tau}(Y)=\rho_{\sigma,\tau}(Z) & \geq & E_Q(-Z|{\cal F}_{\sigma})-\alpha^m_{\sigma,\tau}(Q)\;\;Q\;a.s. \nonumber
\end{eqnarray}
Thus for every $Q \in {\cal M}^1_{\nu,\tau}$, using as before the fact that $E_Q(\alpha^m_{\sigma,\tau}(Q)|{\cal F}_{\nu})$ is finite $Q$ a.s.,
$$\rho_{\nu,\tau}(Y)+\alpha^m_{\nu,\tau}(Q) \geq  E_Q(-Z|{\cal F}_{\nu})\;\;Q\;a.s.$$ 
And then 
$$\rho_{\nu,\tau}(Y) \geq \rho_{\nu,\tau}(Z)$$
Exchanging the roles of $Y$ and $Z$ we get the equality
$\rho_{\nu,\tau}(Y)=\rho_{\nu,\tau}(Z)$.
 From  translation invariance,  $\rho_{\sigma,\tau}(X)=\rho_{\sigma,\tau}(-\rho_{\sigma,\tau}(X)+\rho_{\sigma,\tau}(0))$. We then apply the preceding result  to $Y=X$ and $Z=-\rho_{\sigma,\tau}(X)+\rho_{\sigma,\tau}(0)$. And we get $\rho_{\nu,\tau}(X)=\rho_{\nu,\tau}(-\rho_{\sigma,\tau}(X)+\rho_{\sigma,\tau}(0))$. From hypothesis 2) we get the result.\\

\section{Cadlag modification of a time consistent dynamic risk process}
\label{sec:cadlag} 
Delbaen has proved  (theorem 5.1. of~\cite{D}), that any coherent dynamic risk measure  continuous from above such that $\alpha^m(P)=0$  has a cadlag modification. 
We generalize here this result. We prove that  for every time consistent non degenerate normalized dynamic risk measure  continuous from above  $(\rho_{\sigma,\tau})_{0 \leq \sigma \leq \tau}$, for every $X$ ${\cal F}_{\tau}$-measurable, the process $(\rho_{\sigma,\tau}(X))_{\sigma}$   has a 
modification with cadlag trajectories. We say that the  dynamic risk measure is non degenerate if for any ${\cal F}$-measurable set $A$, $\rho_{0,\infty}(\lambda 1_A)=0$ for every $\lambda>0$ only if $P(A)=0$. 
The aim of the proof is the same as that of Delbaen~\cite{D}. It is related to the construction of the Snell enveloppe. In this proof 
we will make use of the non negativity of the penalty for any normalized dynamic risk measure.

 Denote ${\cal M}^0_{0,T}=\{Q \ll P \;|\; \alpha^m_{0,T}=0\}$
\begin{lemma}
Let $\rho_{\sigma,\tau}$ a time consistent normalized dynamic risk measure continuous from above. Let $T$  a stopping time.  Assume either continuity from below or that ${\cal M}^0_{0,T} \neq \emptyset$.
Let $X \in L^{\infty}({\cal F}_T)$. Let $0 \leq \sigma \leq T$.

Then the process $(\rho_{\sigma,T}(X))_{\sigma \leq T}$ is a supermartingale in the following sense:
$$\forall Q \in {\cal M}^0_{0,T}\;\forall 0 \leq \sigma \leq \tau \leq T\;\;\rho_{\sigma,T}(X)
\geq E_{Q}(\rho_{\tau,T}(X)|{\cal F}_{\sigma})\; Q\;a.s.$$  
\label{lemma1}
\end{lemma}
{\bf Proof}: 
Let $Q  \in {\cal M}^0_{0,T}$. From time consistency, and lemma \ref{lemma0}, 

$\rho_{\sigma, T}(X)=\rho_{\sigma, \tau}(-\rho_{\tau,T}(X))
\geq  (E_Q(\rho_{\tau,T}(X)|{\cal F}_{\sigma})-\alpha^m_{\sigma,\tau}(Q))\;\;Q\;a.s.$.

As the dynamic risk measure is normalized the penalty is always non negative and from the cocycle condition it follows that   $\forall Q \in {\cal M}^0_{0,T}$, $\alpha^m_{\sigma,\tau}(Q)= 0 \;Q\;a.s.$  
$\forall 0 \leq \sigma \leq \tau \leq T $. So $\rho_{\sigma,T}(X)
\geq E_{Q}(\rho_{\tau,T}(X)|{\cal F}_{\sigma})\; \;Q\;a.s.\;\;\forall Q \in {\cal M}^0_{0,T}$.
 
\begin{lemma}
let $(\rho_{\sigma,\tau})$ a time consistent normalized dynamic risk measure continuous from above. Let $T$ a stopping time. Assume that ${\cal M}^0_{0,T}$ contains a probability measure $Q$ equivalent to $P$.  Consider a decreasing sequence of finite stopping times $\sigma_n \leq T$ converging to $\sigma$.\\ 
Then  $E_{Q}(\rho_{\sigma_n,T}(X))$ converges to 
$E_{Q}(\rho_{\sigma,T}(X))$, and $\rho_{\sigma_n,T}(X)$ tends to 
$\rho_{\sigma,T}(X)$ in $L^1(\Omega, {\cal F}_{\infty},P)$.
\label{lemma2}
\end{lemma}

{\bf Proof}: 
Let $X \in L^{\infty}({\cal F}_{T})$. Let $Q \in {\cal M}^0_{0,T}$.
From lemma \ref{lemma0} part ii), 
there is a  sequence $Q_k \in {\cal M}^1_{\sigma,T}(Q)$ such that $E_{Q_k}(-X|{\cal F}_{\sigma})-\alpha^m_{\sigma,T}(Q_k)$ is increasing \;Q \; a.s. and tends to $\rho_{\sigma,T}(X)\; Q\;a.s.$
 
Applying Lemma\ref{lemmaS} 2), as $E_Q(\alpha^m_{\sigma,T}(Q_1)<\infty$, it follows that 
$E_{Q}(E_{Q_k}(-X|{\cal F}_{\sigma})-\alpha_{\sigma,T}(Q_k))$ increases to 
$E_{Q}(\rho_{\sigma,T}(X))$. 

Let $\epsilon >0$.  There is $k_0$ such that for $k>k_0$, $$\epsilon + E_{Q}(\rho_{\sigma,T}(X)) \leq E_{Q}(E_{Q_k}(-X|{\cal F}_{\sigma})-\alpha^m_{\sigma,T}(Q_k))$$
$$=E_{Q}(E_{Q_k}(E_{Q_k}(-X|{\cal F}_{{\sigma}_n})|{\cal F}_{\sigma})-E_{Q}[(\alpha^m_{\sigma,\sigma_n}(Q_k))+E_{Q_k}(\alpha^m_{\sigma_n,T}(Q_k)|{\cal F}_{\sigma})]$$
(applying the cocycle condition of $\alpha^m$).

Now $(E_{Q_k}(-X|{\cal F}_{{\sigma}_n})-\alpha^m_{\sigma_n,T}(Q_k) \leq \rho_{\sigma_n,T}(X)\;\;Q_k\;a.s.$, 
so we obtain
$$\epsilon + E_{Q}(\rho_{\sigma,T}(X)) \leq E_{Q}(E_{Q_k}(\rho_{\sigma_n,T}(X)|{\cal F}_{\sigma}))-E_{Q}(\alpha^m_{\sigma,\sigma_n}(Q_k))$$
$$\leq  E_{Q}(\rho_{\sigma_n,T}(X))+  E_{Q}(\rho_{\sigma_n,T}(X)(E_Q(\frac{dQ_k}{dQ}|{\cal F}_{\sigma_n}))-1))-E_{Q}(\alpha^m_{\sigma,\sigma_n}(Q_k))$$
As the restriction of $Q_k$ to ${\cal F}_{\sigma}$ is equal to $Q$, 
for $k$ fixed,\\
 ${E_{Q}(\frac{dQ_k}{dQ}|{\cal F}_{\sigma_n})}\rightarrow 1$ in $ L^1$, as $n \rightarrow  \infty$, and $||\rho_{\sigma_n,T}(X)||_{\infty}\leq ||X||_{\infty}$, 
so from the dominated convergence theorem, 
$E_{Q}(\rho_{\sigma_n,T}(X)({E_{Q}(\frac{dQ_k}{dQ}|{\cal F}_{\sigma_n})})-1) \rightarrow 0$
as $n \rightarrow \infty$. 

Furthermore $E_Q( \alpha^m_{\sigma,\sigma_n}(Q_k)) \geq 0$ and it follows from the preceding lemma that 
$E_{Q}(\rho_{\sigma,T}(X)) \geq E_{Q}(\rho_{\sigma_n,T}(X))$ for every $n$. So we have the equality
$$E_{Q}(\rho_{\sigma,T}(X)) = \lim_{n \rightarrow \infty} E_{Q}(\rho_{\sigma_n,T}(X))$$ 

We apply the modification theorem (theorem 4 p. 76 in~\cite{DM}) to the $Q$-supermartingale $\rho_{\sigma,T}(X),(\rho_{\sigma_n,T}(X))_n$. 
This $Q$-supermartingale has a\\ 
 $(\Omega,{\cal F}_{\infty},{\cal F}_t,Q)$-modification with cadlag trajectories. And thus applying the dominated convergence 
theorem to this modification 
(as $||\rho_{\sigma_n,T}(X)||_{\infty}\leq ||X||_{\infty})$, $\rho_{\sigma_n,T}(X)$ tends to 
$\rho_{\sigma,T}(X)$ in $L^1(\Omega, {\cal F}_{\infty},Q)=L^1(\Omega, {\cal F}_{\infty},P)$.

\begin{theorem}
Let $(\rho_{\sigma,\tau})_{0 \leq \sigma \leq \tau}$ be a normalized dynamic 
risk measure continuous from above. Let $T$ a stopping time. Assume that there is in ${\cal M}^0_{0,T}$ 
a probability measure $Q$ equivalent to $P$. 
Let $X \in L^{\infty}({\cal F}_T)$.

Then  there is a cadlag $Q$-supermartingale process $Y$ 
such that for every finite stopping time $\sigma \leq \tau$,  
$\rho_{\sigma,T}(X)=Y_{\sigma}$ in $L^{\infty}(\Omega,{\cal F}_{\sigma},P)$
\label{thmC}
\end{theorem}
{\bf Proof}:
The proof of this theorem is similar to that of lemma 5.8. of~\cite{D}, it is based on the 
modification theorem (theorem 2 p. 73 of~\cite{DM}).\\
Let $Q$ equivalent to $P$ $Q\in {\cal M}^0_{0,T}$; Let $X \in  L^{\infty}(\Omega,{\cal F}_T,P)$.   From the preceding lemmas,  
$(\rho_{\sigma,T}(X))_{\sigma <T}$ is a supermartingale for $Q$ 
and $\sigma \rightarrow \rho_{\sigma,T}(X)$ is right continuous in $L^1$.\\
Apply the modification theorem to  the set of  rational numbers.
This gives a cadlag $Q$-supermartingale process $Y$ such that 
$\forall t \in \Q^+\;\;Y_t= \rho_{t,\infty}(X)\;\;Q\;a.s.$\\
Let $\sigma$ a finite stopping time. $\sigma$ is the decreasing limit of a sequence $\sigma_n$ of finite stopping times with rational values.
$Y_{\sigma_n}= \rho_{\sigma_n,\infty}(X)\;\;Q\;a.s.$.\\
Taking the limit in $L^1(\Omega,{\cal F}_{\infty},Q)=L^1(\Omega,{\cal F}_{\infty},P)$, 
applying the lemma \ref{lemma2}, we get the result.

 Now we want to find sufficient conditions on the dynamic risk measure in order to insure that the hypothesis of theorem \ref{thmC} are satisfied.

Recall  the following definition introduced by Peng~\cite{Peng} for non linear expectations:

\begin{definition}
A dynamic risk measure is  strictly monotone if  $\forall (Y,Z) \in L^{\infty}(\Omega,{\cal F},P)$, $Y \geq Z$ and $\rho_{0,\infty}( Y)=\rho_{0,\infty}( Z)$  implies $Y=Z$.
\label{def}
\end{definition}
We introduce a weaker notion that we call non degenerate.
\begin{definition}
A dynamic risk measure is non degenerate if  $\forall A \in {\cal F}_{\infty}$,  
$\rho_{0,\infty}(\lambda 1_A)=\rho_{0,\infty}(0)\;\;\forall \lambda \in {\R_{+}}^*$ implies $P(A)=0$.
\label{defbis}
\end{definition}
\begin{lemma}
Let  $(\rho_{\sigma,\tau})_{0 \leq \sigma \leq \tau}$ 
be a non degenerate normalized time consistent dynamic risk measure continuous from above. 
Let $\tau$  a stopping time. Every probability measure in   
${\cal M}^0_{0,\tau}$ is equivalent to $P$ on $(\Omega,{\cal F}_{\tau})$.
\label{lemma3.4}
\end{lemma}
{\bf Proof}:
Let $A \in {\cal F}_{\tau}$ such that $P(A)>0$  
As $\rho$ is non degenerate,  
there  is  $\lambda \in {\R_{+}}^* $ such that $\rho _{0,{\infty}}(\lambda 1_A)< 0$. From normalization and time consistency $\rho _{0,{\tau}}(\lambda 1_A)
=\rho _{0,{\infty}}(\lambda 1_A)$. Then  
$\forall Q \in {\cal M}^0_{0,{\tau}}$, 
$\rho_{0,{\tau}}(\lambda 1_A) \geq -\lambda Q(A)$.
So $Q(A)>0$. As $Q\ll P$, it follows that $Q$ is equivalent to $P$.

Q.e.d.
\begin{corollary}

  Consider a time consistent dynamic risk measure continuous from above, normalized and non degenerate. Assume that ${\cal M}^0_{0,\infty} \neq \emptyset$ or that the risk measure is continuous from below.  Let $X \in L^{\infty}({\cal F})$. For every $Q \in {\cal M}^0_{0,\infty}$, there is a cadlag $Q$-supermartingale process $Y$ 
such that for every finite stopping time  $\sigma$,  
$\rho_{\sigma,\infty}(X)=Y_{\sigma}$ in $L^{\infty}(\Omega,{\cal F}_{\sigma},P)$ 
\label{cor}
\end{corollary}
 In the case of continuity from below, from equation \ref{eq_1}, as $\rho_{0,\infty}(0)=0$, we obtain that ${\cal M}^0_{0,\infty}(0)=0$.
 The corollary  is then an immediate consequence of theorem \ref{thmC} and of lemma \ref{lemma3.4}.

\section{Examples}

 We give here some examples of time consistent dynamic risk measures.
\subsection{Dynamic risk measures defined {\it ex ante}}
\subsubsection*{Entropic dynamic risk measure with threshold}

Consider $(g_{s,t})_{0 \leq s <t}$ a  family of stictly positive bounded ${\cal F}_s$-measurable functions such that $\ln(g_{s,t})$ is essentially bounded.
Consider the entropic dynamic risk measure defined as follows:

Let $0\leq s \leq t$. For every $X \in  L^{\infty}({\cal F}_t)$
\begin{eqnarray}
\rho_{s,t}(X) & = & \essinf\{Y \in {\cal E}_{{\cal F}_s}\;/\;E(e^{[-\alpha (X+Y)]}|
{\cal F}_s) \leq g_{s,t}\} \nonumber \\
 & = & \frac{1}{\alpha} [\ln E(e^{-\alpha X}|{\cal F}_s)-\ln(g_{s,t}) ]. \nonumber
\end{eqnarray}

$\rho_{s,t}$ is normalized iff $g_{s,t}=1 \;\forall s\leq t$.

From~\cite{BN01}, the minimal penalty is  $\alpha^m_{s,t}(Q)=\frac{1}{\alpha}(E_P(\ln(\frac{dQ}{dP})\frac{dQ}{dP}|{\cal F}_s)-\ln(g_{s,t}))$.

The entropic dynamic risk measure is time consistent if and only if  the functions $g_{s,t}$ are ${\cal F}_0$-measurable  i.e. a.s. constant and satisfy the relation 
$\forall \; r,s,t; \; 0 \leq r  \leq s \leq t, \;\; \ln(g_{r,t})=\ln(g_{r,s})+\ln(g_{s,t})\;a.s.$ (see \cite{BN02}). In particular if we assume that there is a strictly positive real valued continuous function $h$ such that 
$\forall (s,t)\; g_{s,t}=h(t-s)$ then the associated dynamic risk measure is time-consistent if and only if there is a real number $\lambda$ such that 
$g_{s,t}=e^{\lambda(t-s)}$.
$\rho_{s,t}$ is normalized iff $g_{s,t}=1 \;\forall s\leq t$.\\
The time consistency of the usual entropic dynamic risk measure ($g_{s,t}=1$) has been studied in  ~\cite{BEK}  ~\cite{DS},~\cite{KS}. With thresholds it has been studied in ~\cite{BN02} 
\subsubsection*{Examples based on BSDE}
The dynamic risk measures coming from B.S.D.E. are time consistent (see \cite{Peng}, \cite{RG}, \cite{BEK} and \cite{KS}).
 Consider $(\Omega,{\cal F},{\cal F}_t,P)$ where ${\cal F}_t$ is the augmented filtration of a $d$ dimensional Brownian motion. Assume that the driver $g(t,z)$ satisfies $g(t,0)=0$ and is  continuous and convex (in $z$),  and satisfies the condition of quadratic growth. The associated BSDE, 
\begin{eqnarray}
-dY_t & = & g(t,Z_t)dt -Z_t^*dB_t \nonumber \\
Y_T & = &  X \nonumber 
\end{eqnarray}
has a solution which gives rise to a  dynamic risk measure $\rho_{s,T}(-X)=Y_s$.
Barrieu and El Karoui~\cite{BEK}, 
section~7.3, have computed the minimal penalty  associated to this dynamic risk measure.

\subsection{Examples constructed from a stable set of probability measures}
 Recall the definition of  stability for a set  of probability measures (cf ~\cite{D})
\begin{definition}
A set ${\cal Q}$ of probability measures all equivalent to $P$ is stable
 if for every stopping times, $\nu \leq \sigma \leq \tau$,
For every $Q \in  {\cal Q}$, for every  $R \in  {\cal Q}$, there is $S \in  {\cal Q}$ such that 
$$\forall f \in L^{\infty}({\cal F}_{\tau}), \;E_S(f|{\cal F}_{\nu})=E_R(E_Q(f|{\cal F}_{\sigma})|{\cal F}_{\nu})\;P.a.s.$$
\end{definition}
Recall the definition of locality for a penalty function
\begin{definition}
 A  penalty function $\alpha_{\sigma,\tau}$   
is local if for every stopping times $\sigma \leq \tau$, for every $A$ 
${\cal F}_{\sigma}$-measurable, if $E_{Q_1}((X1_A|{\cal F}_{\sigma})=E_{Q_2}(X1_A|{\cal F}_{\sigma})\;P.a.s. \; \forall X \in L^{\infty}(\Omega,{\cal F}_{\tau},P)$, then $\;1_A \alpha_{\sigma,\tau}(Q_1)=1_A \alpha_{\sigma,\tau}
(Q_2)\;P.a.s.$
\end{definition}
Recall that from~\cite{BN02} theorem 6, any stable family of probability measures and any  local penalty satisfying the cocycle condition lead to a time consitent dynamic risk measure. This is an important way to construct time consistent dynamic risk measures.
\subsubsection*{Examples constructed from continuous martingales}

 $$\;\;{\cal Q}_1=\{Q_M\;;\;\frac{dQ_M}{dP}={\cal E}(M)\;|\; M \;continuous\;P\;  martingale\;;$$
$$ [M,M]_{\infty} \in L^{\infty}(\Omega,{\cal F},P\}$$ 
 is a stable set of probability measures all equivalent to $P$.
$${\cal Q}_2=\{Q_M\;;\;\frac{dQ_M}{dP}={\cal E}(M)\;|\; M \;continuous\;P \;  martingale\;;$$
$$ ||M||_{BMO}^2=sup_{S } ||E_{P}([M,M]_{\infty}-[M,M]_{S^-}|{\cal F}_S)||_{\infty} <\infty \}$$ 
is stable. This follows from lemma 5 of~\cite{BN02}. It is the set of continuous BMO martingales. For the continuous BMO martingales we refer to ~\cite{K}. 
BMO continuous martingales are also used in ~\cite{BEK} in order to compute the penalty associated to a BSDE.\\
Define on ${\cal Q}_i$ (whith $i$ equal to ${1}$ or ${2}$ the penalty function $\alpha$ as follows:
$$\forall 0 \leq \sigma \leq \tau \;\;\alpha_{\sigma,\tau}(Q_M)= E_{Q_M}( \int\limits_{\sigma}^{\tau} b_u d[M,M]_u|{\cal F}_{\sigma})$$
where $b_u$ is a  non negative predictable process.
Then
$$\rho_{\sigma,\tau}(X)=esssup_{ M \in {\cal M}} (E_{Q_M}(-X|{\cal F}_{\sigma})- \alpha_{\sigma,\tau}(Q_M))$$
defines a time-consistent normalized dynamic risk measure on the filtered probability space 
$(\Omega,{\cal F},P,({\cal F}_t)_{0 \leq t \leq \infty})$ 'cf ~\cite{BN02}).  The penalty function associated to $P=Q_0$ is always equal to $0$. Thus for every $X \in {\cal F}_{\tau}$, the process $(\rho_{\sigma,\tau}(X))_{\sigma}$ has a cadlag modification.

\subsubsection*{Examples based on BMO right continuous martingales}
We can generalize the preceeding example considering a stable set ${\cal M}$ of BMO right continuous  martingales  of norm BMO uniformly bounded by $m$. 
For the general theory of BMO martingales we refer to ~\cite{DDM01} and ~\cite{DDM02}.
For example:  
let $M^1,...,M^j$  a family of strongly orthogonal square integrable right continuous martingales in $(\Omega,{\cal F},({\cal F}_t)_{0 \leq t},P)$. Consider $(\Phi_i)_{1 \leq i \leq j}$  non negative predictable processes such that $\Phi_i M^i$ is a BMO martingale of BMO norm  $m^i$. $${\cal M}=\{\sum_{1 \leq i \leq j} H_i.N_i \;|\;H_i \;predictable\;|H_i| \leq \phi_i \;\; a.s.\}$$ is a  set of square integrable BMO martingales  with norm BMO uniformly bounded by $(\sum _{1 \leq i \leq j} (m^i)^2)^\frac{1}{2}=m$.

If $M^1,...,M^j$ are continuous or if $m<\frac{1}{16}$,   ${\cal Q}({\cal M})$ the corresponding set  of probability measures $(Q_M)_{M \in {\cal M}}$  of Radon Nikodym derivative $\frac{dQ_M}{dP}={\cal E}(M)$.
is stable.  Let $b_s$ be a bounded  predictable process. 
Define on  ${\cal Q}({\cal M})$ the penalty function $\alpha$ by:

for every stopping times $0 \leq \sigma \leq \tau$
$$\alpha_{\sigma,\tau}(Q_M)=E_{Q_M}( \int \limits_{\sigma}^{\tau} b_u d[M,M]_u)|{\cal F}_{\sigma})$$
Then from ~\cite{BN02}, $$\rho_{\sigma,\tau}(X)=\esssup_{ Q_M \in {\cal Q}({\cal M})} (E_{Q_M}(-X|{\cal F}_{\sigma})+ \alpha_{\sigma,\tau}(Q_M))$$
defines a time consistent dynamic risk measure.
 If $0 \in {\cal M}$ and $b$ is non negative, the dynamic risk measure is normalized and satisfies $\alpha_{\sigma,\tau}(P)=0$. 
 Variants of this example con be found in ~\cite{BN02}. For example $b_u$ can depend on the $H_i$.\\
We have thus constructed a large class of examples of time consistent dynamic risk processes with cadlag but not necessary continuous trajectories. These examples generalize the BSDE. Indeed it follows from ~\cite{BEK} that the minimal penalty associated to a BSDE with a strictly convex driver is of that form (with continuous BMO martingales.

\section{Conclusion}

 We have extended to the  general case of dynamic risk measures continuous from above the characterization of time consistency in terms of a cocycle condition\\
 We have also extended the supermartingale property. \\ 
We have proved that for any  non degenerate normalized time consistent dynamic risk
measure $\rho_{\sigma,\tau}$  continuous from  above , under the condition of the existence  of  a probability measure $Q$ equivalent to $P$ of zero penalty , for every bounded measurable $X$, there is a cadlag $Q$-supermartingale  process $Y$ such that for every stopping time $\sigma$, $Y_{\sigma}= \rho_{\sigma,{\infty}}(X)$ $P$\; a.s. When the dynamic risk measure is continuous from below, this condition  is always satisfied.\\
 Any stable set of equivalent probability measures and any local penalty satisfying the cocycle condition give rise to a time consistent dynamic risk measure. Using the theory of right continuous BMO martingales we then construct examples of time consistent dynamic risk processes with possible jumps.
The application of dynamic risk measuring to pricing in incomplete markets will be the subject of a forthcoming paper \cite{BN04}.

\subsection*{Acknowledgements}
I thank Nicole El Karoui and Rama Cont for useful comments, and Nicole El Karoui
for a critical reading of the manuscript.

\end{document}